\documentclass[12pt]{article}
\usepackage{amsfonts}

\usepackage{latexsym}
\usepackage{amssymb}
\usepackage{epsfig}
\usepackage{amsmath}
\usepackage{amsthm}
\usepackage[mathscr]{eucal}
\usepackage{subfigure}
\usepackage{graphicx}

\newtheorem{Lemma}{Lemma}

\newtheorem{Theorem}[Lemma]{Theorem}

\newtheorem{Definition}{Definition}

\renewcommand{\qed}{\hfill{\ \ \rule{2mm}{2mm}} \vspace{0.2in}}

\begin{document}

\title{Outermost boundaries for star-connected components in percolation}
\author{ \textbf{Ghurumuruhan Ganesan}
\thanks{E-Mail: \texttt{gganesan82@gmail.com} } \\
\ \\
NISER, Bhubaneshwar, India}
\date{}
\maketitle

\begin{abstract}



Tile \(\mathbb{R}^2\) into disjoint unit squares \(\{S_k\}_{k \geq 0}\) with the origin being the centre of \(S_0\)
and say that \(S_i\) and \(S_j\) are star-adjacent if they share a corner and plus-adjacent if they share an edge.
Every square is either vacant or occupied. If the occupied
plus-connected component \(C^+(0)\) containing the origin is finite, it is known
that the outermost boundary \(\partial^+_0\) of \(C^+(0)\) is a unique cycle surrounding the origin.
For the finite occupied star-connected component \(C(0)\) containing the origin,
we prove in this paper that the outermost boundary \(\partial_0\) is a unique connected graph
consisting of a union of cycles \(\cup_{1 \leq i \leq n} C_i\) with mutually disjoint
interiors. Moreover, we have that each pair of cycles in \(\partial_0\) share at
most one vertex in common and we provide an inductive procedure to obtain a circuit containing
all the edges  of \(\cup_{1 \leq i \leq n} C_i.\) This has applications for contour analysis
of star-connected components in percolation.

\vspace{0.1in} \noindent \textbf{Key words:} Star connected components, outermost boundary, union of cycles.

\vspace{0.1in} \noindent \textbf{AMS 2000 Subject Classification:} Primary:
60J10, 60K35; Secondary: 60C05, 62E10, 90B15, 91D30.
\end{abstract}

\bigskip

\renewcommand{\theequation}{\thesection.\arabic{equation}}
\setcounter{equation}{0}
\section{Introduction} \label{intro}

Tile \(\mathbb{R}^2\) into disjoint unit squares \(\{S_k\}_{k \geq 0}\) with origin being the centre of~\(S_0.\) We say \(S_1\) and \(S_2\) are \emph{adjacent} or \emph{star-adjacent} if they share a corner between them. We say that squares \(S_1\) and \(S_2\) are \emph{plus-adjacent}, if they share an edge between them. Here we follow the notation of Penrose (2003). Suppose every square is assigned one of the two states: occupied or vacant. In many applications like for example, percolation, it is of interest to determine the outermost boundary of the plus-connected or star-connected components containing the origin. We make formal definitions below. The case of plus-connected components is well studied (Bollobas and Riordan (2006), Penrose (2003)) and in this case, the outermost boundary is simply a cycle containing the origin. Our main result is that the outermost boundary for the star-connected component is a connected union of cycles with disjoint interiors. 

Let \(C(0)\) denote the star-connected occupied component containing the origin and throughout we assume that \(C(0)\) is finite. Thus if \(S_0\) is vacant then \(C(0) = \emptyset.\) Else \(S_0 \in C(0)\) and if \(S_1, S_2 \in C(0)\) there exists a sequence of distinct occupied squares \((Y_1,Y_2,...,Y_t)\) all belonging to \(C(0),\) such that \(Y_i\) is adjacent to \(Y_{i+1}\) for all \(i\) and \(Y_1 = S_1\) and \(Y_t = S_2.\) Let \(G_C\) be the graph with vertex set being the set of all corners of the squares \(\{S_k\}_k\) in \(C(0)\) and edge set consisting of the edges of the squares \(\{S_k\}_k\) in~\(C(0).\)

Two vertices \(u\) and \(v\) are said to be adjacent in \(G_C\) if they share an edge between them. We say that an edge \(e\) in \(G_C\) is adjacent to square \(S_k\) if it is one of the edges of \(S_k.\) We say that \(e\) is a \emph{boundary edge} if it is adjacent to a vacant square and is also adjacent to an occupied square. A \emph{path} \(P\) in \(G_C\) is a sequence of distinct vertices \((u_0,u_1,...,u_t)\) such that \(u_i\) and \(u_{i+1}\) are adjacent for every \(i.\)  A \emph{cycle} \(C\) in \(G_C\) is a sequence of distinct vertices \((v_0,v_1,...,v_m,v_0)\) starting and ending at the same point such that \(v_i\) is adjacent to \(v_{i+1}\) for all \(0 \leq i \leq m-1\) and \(v_m\) is adjacent to \(v_0.\) A \emph{circuit} \(C'\) in \(G_C\) is a sequence of vertices \((w_0,w_1,...,w_r,w_0)\) starting and ending at the same point such that \(w_i\) is adjacent to \(w_{i+1}\) for all \(0 \leq i \leq r-1,\) \(w_r\) is adjacent to \(w_0\) and no edge is repeated in \(C'.\) Thus vertices may be repeated in circuits and for more related definitions, we refer to Chapter~1, Bollobas (2001).

Any cycle \(C\) divides the plane \(\mathbb{R}^2\) into two disjoint connected regions. As in Bollobas and Riordan (2006), we denote the bounded region to be the \emph{interior} of \(C\) and the unbounded region to be the \emph{exterior} of \(C.\) We have the following definition.
\begin{Definition} We say that edge \(e\) in \(G_C\) is an \emph{outermost boundary} edge of the component \(C(0)\) if the following holds true for every cycle \(C\) in \(G_C:\) either \(e\) is an edge in \(C\) or \(e\) belongs to the exterior of \(C.\)

We define the outermost boundary \(\partial _0\) of \(C(0)\) to be the set of all outermost boundary edges of \(G_C.\)
\end{Definition}
Thus outermost boundary edges cannot be contained in the interior of any cycle in \(G_C.\) Our main result is the following.
\begin{Theorem}\label{thm3} Suppose \(C(0)\) is finite.
The outermost boundary \(\partial_0\) of \(C(0)\) is a unique set of cycles \(C_1,C_2,...,C_n\) in \(G_C\) with the following properties:\\
(i) The graph \(\cup_{1 \leq i \leq n}C_i\) is a connected subgraph of \(G_C.\)\\
(ii) If \(i \neq j,\) the cycles \(C_i\) and \(C_j\) have disjoint interiors and share at most one vertex.\\
(iii) Every square \(S_k \in C(0)\) is contained in the interior of some cycle \(C_{j}.\)\\
(iv) If \(e \in C_{j}\) for some \(j,\) then \(e\) is a boundary edge of \(C(0)\) adjacent to an occupied square of \(C(0)\) in the
interior of \(C_j\) and also adjacent to a vacant square in the exterior.

Moreover, there exists a circuit \(C_{out}\) containing every edge of \(\cup_{1 \leq i \leq n} C_i.\)
\end{Theorem}
The outermost boundary \(\partial_0\) is therefore also an Eulerian graph with \(C_{out}\) denoting the corresponding Eulerian circuit (for definitions, we refer to Chapter 1, Bollobas (2001)). We remark that the above result also provides a more detailed justification of the statement made about the outermost boundary and the corresponding circuit in the proof of Lemma~3 of Ganesan (2013). Using the above result, we also obtain the outermost circuit that is used to construct the top-down crossing in oriented percolation in a rectangle in Ganesan (2015). 


The proof of the above result also obtains the outermost boundary cycle in the case of plus-connected components. We recall that \(S_1\) and \(S_2\) are \emph{plus-adjacent} if they share an edge between them. Analogous to the star-connected case, we define \(C^+(0)\) to be the plus-connected component containing the origin and define the graph \(G^+_C\) consisting of edges and corners of squares in \(C^+(0).\) We have the following.
\begin{Theorem}\label{thm2} Suppose \(C^+(0)\) is finite. The outermost boundary \(\partial^+_0\) of \(C^+(0)\) is unique cycle \(C^+_{out}\) in \(G^+_C\) with the following property:\\
(i) All squares of \(C^+(0)\) are contained in the interior of \(C^+_{out}.\)\\
(ii) Every edge in \(C^+_{out}\) is a boundary edge adjacent to an occupied square of \(C^+(0)\) in the interior of \(C^+_{out}\) and a vacant square in the exterior.
\end{Theorem}
This is in contrast to star-connected components which may contain multiple cycles in the outermost boundary.  

To prove Theorem~\ref{thm3}, we use the following intuitive result about merging cycles. Analogous to \(G_C,\) let \(G\) be the graph with vertex set being the corners of the squares \(\{S_k\}_k\) and edge set being the edges of the squares \(\{S_k\}_k.\)
\begin{Theorem}\label{thm1} Let \(C_1\) and \(C_2\) be cycles in \(G\) that have more than one vertex in common. There exists a unique cycle \(C_3\) consisting only of edges of \(C_1\) and \(C_2\) with the following properties:\\
(i) the interior of \(C_3\) contains the interior of both \(C_1\) and \(C_2,\)\\
(ii) if an edge \(e\) belongs to \(C_1\) or \(C_2,\) then either \(e\) belongs to \(C_3\) or is contained in its interior.

Moreover, if \(C_2\) contains at least one edge in the exterior of \(C_1,\) then the cycle \(C_3\) also contains an edge of \(C_2\) that lies in the exterior of \(C_1.\)
\end{Theorem}
The above result essentially says that if two cycles intersect at more than one point, there is a innermost cycle containing both of them in its interior. We provide an iterative construction for obtaining the cycle \(C_3,\) analogous to Kesten~(1980) for crossings, in Section~\ref{pf1}.

The paper is organized as follows: In Section~\ref{pf2}, we prove Theorem~\ref{thm3} and in Section~\ref{pf1}, we prove Theorem~\ref{thm2} and Theorem~\ref{thm1}.


\section{Proof of Theorem~\ref{thm3}} \label{pf2}
\emph{Proof of Theorem~\ref{thm3}}: The first step is to obtain large cycles surrounding each occupied square in \(C(0).\) We have the following Lemma.

\begin{Lemma}\label{outer} For every \(S_k \in C(0),\) there exists a unique cycle \(D_k\) satisfying the following properties:\\
(a) \(S_k\) is contained in the interior of \(D_k,\)\\
(b) every edge in the cycle \(D_k\) is a boundary edge adjacent to one occupied square of \(C(0)\) in the interior and one vacant square in the exterior and\\
(c) if \(C\) is any cycle in \(G_C\) that contains \(S_k\) in the interior, then every edge in \(C\) either belongs to \(D_k\) or is contained in the interior.
\end{Lemma}
We denote \(D_k\) to be the outermost boundary cycle containing the\\square \(S_k.\) We prove all statements at the end.

We claim that the set of distinct cycles in the set \({\cal D} := \cup_{S_k \in C(0)} \{D_k\}\) is the desired outermost boundary \(\partial_0\) and satisfies the conditions (i)-(iv) mentioned in the statement of the theorem. By construction, we have that (iii) and (iv) are satisfied. To see that (ii) holds, we suppose that \(D_{k_1} \neq D_{k_2}\) and that \(D_{k_1}\) and \(D_{k_2}\) meet at more than one vertex. We know that \(D_{k_2}\) is not completely contained in \(D_{k_1}.\) Thus \(D_{k_2}\) contains at least one edge in the exterior of \(D_{k_1}.\)  From Theorem~\ref{thm1}, we obtain a cycle \(D'_{12}\) containing both \(D_{k_1}\) and \(D_{k_2}\) in the interior and containing an edge \(e\) present in \(D_{k_2}\) but not in \(D_{k_1}\) or its interior. The cycle \(D'_{12}\) satisfies condition (a) in Lemma~\ref{outer} above and thus contradicts the assumption that \(D_{k_1}\) satisfies~(c). Thus \(D_{k_1}\) and \(D_{k_2}\) cannot meet at more than one vertex.

Also (i) holds, because of the following reason. First we note that by construction \(G_C\) is connected; let \(u_1\) and \(u_2\) be vertices in \(G_C.\) Each \(u_i, i = 1,2\) is a corner of an occupied square \(S_i \in C(0)\) and by definition, \(S_1\) and \(S_2\) are star-connected via squares in \(C(0).\) Thus there exists a path in \(G_C\) from \(u_1\) to \(u_2.\)

To see that \({\cal D}\) is a connected subgraph of \(G_C,\) we let \(v_1\) and \(v_2\) be vertices in \({\cal D}\) that belong to cycles \(D_{r_1}\) and \(D_{r_2},\) respectively, for some \(r_1\) and~\(r_2.\) If \(r_1 =\)~\(r_2,\) then \(v_1\) and \(v_2\) are connected by a path in \(D_{r_1}= (z_1 = v_1,z_2,...,z_n,z_1).\) If \(r_1 \neq r_2,\) let \(P_{12} = (w_1 = v_1,w_2,...,w_{t-1},w_{t} = v_2)\) be a path from \(v_1\) to \(v_2\) in~\(G_C.\) We iteratively construct a path \(P'_{12}\) from \(P_{12}\) using only edges of cycles in \({\cal D}.\) We first note that since (iii) holds, every edge in \(P_{12}\) either belongs to a cycle in \({\cal D}\) or is contained in the interior of some cycle in \({\cal D}.\) Let \(i_1\) be the first time \(P_{12}\) leaves \(D_{r_1};\) i.e., let \(i_1 = \min\{i \geq 1 : w_{i+1} \text{ belongs to exterior of } D_{r_1}\}.\)

The edge formed by the vertices \(w_{i_1}\) and \(w_{i_1+1}\) belongs to some cycle \(D_{s_1} = (x_1 = w_{i_1},x_2,...,x_r,x_1)\) or is contained in its interior. Since the cycles \(D_{r_1}\) and \(D_{s_1}\) have disjoint interiors, this necessarily means \(D_{s_1}\) and \(D_{r_1}\) meet at \(w_{i_1}.\) Defining \(T_1 = (z_1 = v_1,z_2,...,z_{j_1} = w_{i_1}),\) we note that \(T_1\) is a path consisting only of edges in the cycle \(D_{r_1}\) and containing the vertex \(z_1 = v_1.\) Repeating the same procedure above, we obtain another path \(T_2 = (w_{i_1} = x_1,x_2,...,x_{j_2} = w_{i_2})\) contained in \(D_{s_1},\) where, as before, \(i_2 = \min\{i \geq i_1+1 : w_{i+1} \text{ belongs to exterior of } D_{s_1}\}\) denotes the first time \(P_{12}\) leaves~\(D_{s_1}.\) We continue this procedure for a finite number of steps \(m,\) until we reach~\(v_2.\) By construction, the path \(T_i\) obtained at step~\(i, 2 \leq i \leq m\) is connected to \(\cup_{1 \leq j \leq i-1} T_j.\) The final union of paths \(\cup_{1 \leq i \leq m}T_{i}\) is therefore a connected graph containing only edges in \({\cal D}\) and contains \(v_1\) and \(v_2.\) 



It remains to see that an edge \(e\) belongs to the outermost boundary if and only if it belongs to some cycle in \({\cal D}.\) If \(e\) is an edge in a cycle \(D_k \in {\cal D}\) we have that \(e\) is adjacent to an occupied square \(S_e\) contained in the interior of \(D_k\) and a vacant square \(S'_e\) in the exterior. If there exists a cycle \(C\) in \(G_C\) that contains \(e\) in the interior, we then have that both \(S_e\) and \(S'_e\) are contained in the interior of \(C.\) Since \(S'_e\) is exterior to~\(D_k,\) the cycle \(C\) contains at least one edge in the exterior of \(D_k.\) But if \(D_e\) denotes the outermost cycle containing~\(S_e,\) then by the discussion in the first paragraph, we must have that \(D_e = D_k.\) And thus every edge of \(C\) either belongs to \(D_e\) or is contained in the interior of \(D_e\) which leads to a contradiction.

We also see that no other edge apart from edges of cycles in \({\cal D}\) can belong to the outermost boundary since if \(e_1 \notin {\cal D},\) then \(e_1\) is necessarily contained in the interior of some cycle \(D_r \in {\cal D}.\)

Finally, to obtain the circuit we compute the cycle graph \(H_{cyc}\) as follows: let \(E_1,E_2,...,E_n\) be the distinct outermost boundary cycles in \({\cal D}.\) Represent \(E_i\) by a vertex \(i\) in \(H_{cyc}.\) If \(E_i\) and \(E_j\) share a corner, we draw an edge between \(i\) and \(j.\) We have the following lemma.
\begin{Lemma}~\label{hcyc} We have that the graph \(H_{cyc}\) described above is a tree.
\end{Lemma}
We provide the proof of the above at the end.

We then obtain the circuit via induction on the number of vertices \(n\) of~\(H_{cyc}.\) For \(n = 1,\) it is a single cycle. Suppose we obtain the circuit of all cycle graphs containing at most \(k\) vertices and  let \(H_{cyc}\) be a cycle graph containing \(k+1\) vertices. To obtain the circuit for \(H_{cyc},\) we pick a leaf \(q\) of \(H_{cyc}\) and apply induction assumption on the cycle graph \(H'_{cyc} = H_{cyc}\setminus q.\) To fix a procedure, we choose \(q\) such that the corresponding boundary cycle \(E_q\) contains a square \(S_j\) of least index \(j\) in its interior. We have that \(H'_{cyc}\) is connected and has \(k\) vertices and thus has a circuit  \(C_{k} = (c_1,c_2,...,c_r,c_1)\) containing all edges of every cycle in~\(H'_{cyc}.\) Let \(C_k\) meet the cycle \(E_q = (d_1,d_2,...,d_i,d_1)\) at \(d_t = c_1.\) We then form the new circuit \(C_{k+1} = (d_1,d_2,...,d_t = c_1,c_2,...,c_r,c_1 = d_t,d_{t+1},...,d_i,d_1),\) which contains all edges of every cycle in~\(H_{cyc}.\) \(\qed\)

\emph{Proof of Lemma~\ref{outer}}: We note that if there exists such a \(D_k,\) then it is unique by definition. Let \({\cal E}\) be the set of all cycles in \(G_C\) satisfying condition~(a); i.e., if \(C\) is a cycle containing \(S_k\) in its interior then \(C \in {\cal E}.\) The set \({\cal E}\) is not empty since the cycle formed by the four edges of \(S_k\) belongs to \({\cal E}.\) We merge cycles in \({\cal E}\) two by two using Theorem~\ref{thm1} to obtain the desired cycle~\(D_k.\) We first pick a cycle \(F_1\) in \({\cal E}\) using a fixed procedure; for example, using an analogous iterative procedure as described in Section~1 of Ganesan (2014) for choosing paths. 

We again use the same procedure to pick a cycle \(F_2\) in \({\cal E} \setminus F_1\) and from Theorem~\ref{thm1}, obtain a cycle \(F'_1\) consisting of only edges of \(F_1\) and \(F_2\) and containing both \(F_1\) and~\(F_2\) in its interior. The cycle \(F'_1\) also satisfies (a) and thus belongs to~\({\cal E}.\) Therefore, if \({\cal E}\) has \(t\) cycles, then \({\cal E}_1 := ({\cal E} \setminus \{F_1,F_2\}) \cup F'_1\) has at most \(t-1\) cycles; if \(F_1\) contains an edge in the exterior of \(F_2\) and the cycle \(F_2\) also contains an edge in the exterior of \(F_1,\) then \({\cal E}_1\) has \(t-2\) cycles. Else \(F'_1\) is either \(F_1\) or \(F_2\) and the set \({\cal E}_1\) therefore contains \(t-1\) cycles.

By construction, every cycle in \({\cal E}\) is either a cycle in \({\cal E}_1\) or is contained in the interior of a cycle in~\({\cal E}_1.\) Therefore, if \({\cal E}_1\) contains one cycle, it is the desired outermost boundary cycle \(D_k.\) Else we repeat the above procedure with~\({\cal E}_1\) and obtain another set \({\cal E}_2\) containing at most \(t-2\) cycles and again with the property that every cycle in \({\cal E}\) is either a cycle in \({\cal E}_2\) or is contained in the interior of a cycle in~\({\cal E}_2.\) Continuing this process, we are finally left with a single cycle \(C_{fin}.\) By construction it satisfies (a) and (c). It only remains to see that (b) is true.

Suppose there exists an edge \(e\) of \(C_{fin}\) that is not a boundary edge. Since \(e\) is an edge of \(G_C,\) we then have that \(e\) is adjacent to two occupied squares \(S_1\) and \(S_2,\) with one of the squares, say \(S_1,\) contained in the interior of \(C_{fin}\) and the other square \(S_2,\) contained in the exterior. The cycle \(C_2\) containing the four edges of the square \(S_2\) and the cycle \(C_{fin}\) have the edge \(e\) in common and thus more than one vertex in common. Since \(C_2\) contains at least one edge in the exterior of \(C_{fin},\) we use Theorem~\ref{thm1} to obtain a larger cycle \(C'_2\) containing both \(C_{fin}\) and \(C_2\) in the interior. The cycle \(C'_2\) contains at least one edge not in \(C_{fin}.\) But since \(C_{fin}\) satisfies (c), this is a contradiction. Thus every edge \(e\) of \(C_{fin}\) is a boundary edge.

By the same argument above, we also see that the edge \(e\) cannot be adjacent to an occupied square in the exterior of \(C_{fin}.\) Thus \(e\) is adjacent to an occupied square in the interior and a vacant square in the exterior. \(\qed\)


\emph{Proof of Lemma~\ref{hcyc}}: We already have that \(H_{cyc}\) is connected. It is enough to see that it is acyclic. Before we prove that, we make the following observation. Consider a path \(P = (i_1,i_2,...,i_m)\) in \(H_{cyc}.\) We see that any vertex in \(E_{i_1}\) and any vertex in \(E_{i_m}\) is connected by a path consisting only of edges of the cycles \(\{E_{i_k}\}_{1 \leq k \leq m}.\)


\begin{figure*}
\centerline{\subfigure[]{\includegraphics[width=3.5in, trim= 100 250 100 275, clip=true]{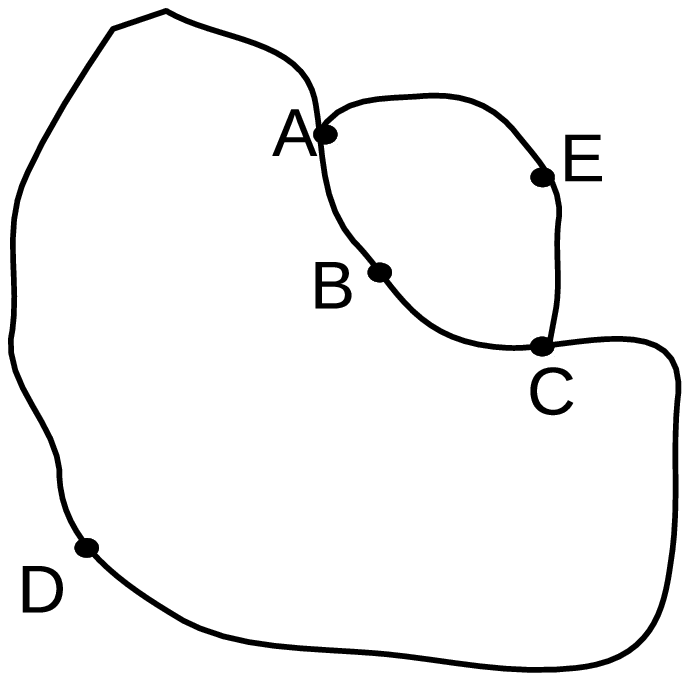}
\label{fig_first_case}}
\hfil
\subfigure[]{\includegraphics[width=3.5in, trim= 100 350 100 200, , clip=true]{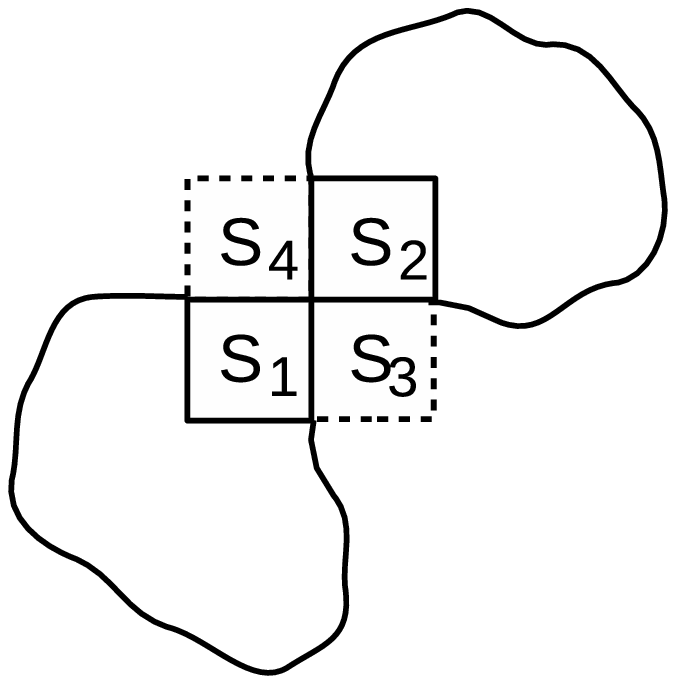}
\label{fig_second_case}}}
\caption{(a) Merging cycle \(ABCDA\) with the segment \(AEC.\) (b) Only two cycles can meet at a single point.}
\label{cyc_fig}
\end{figure*}

Suppose \(H_{cyc}\) contains a cycle \(C = (r_1, r_2,...,r_s, r_1).\) Let the boundary cycle \(E_{r_1}  = (u_1,u_2,...,u_m,u_1)\) meet \(E_{r_2}\) at \(u_{1}\) and \(E_{r_s}\) at \(u_{j}.\) We have that \(j\neq 1\) since three boundary cycles cannot meet at a point. This is illustrated in Figure~\ref{fig_second_case}. The occupied square \(S_1\) belongs to \(E_{r_2}\) and the occupied square \(S_2\) belongs to \(E_{r_s}.\) It is necessary that the squares \(S_3\) and \(S_4\) are vacant and thus cannot be on the boundary of any other cycle.

Let \(P_1\) and \(P'_1\) be the two segments of \(E_{r_1}\) starting at \(u_{1}\) and ending at~\(u_{j}.\) Since \(u_{1} \in E_{r_2}\) and \(u_{j} \in E_{r_s},\) we have by the observation made in the first paragraph that there exists a path \(P_2\) from \(u_1\) to \(u_j,\) consisting only of edges in \(\{E_{r_i}\}_{2 \leq i \leq s}.\) This path necessarily lies in the exterior of \(E_{r_1}\) and is illustrated in Figure~\ref{fig_first_case}. Here \(ABCDA\) represents the cycle \(E_{r_1},\) the path \(P_1\) is the segment \(ADC\) and the path \(P'_1\) is the segment \(ABC.\) The path \(P_2\) is denoted by the exterior segment \(AEC.\)

Thus it is necessary that either the cycle \(C_{12}\) formed by \(P_1 \cup P_2\) contains \(P'_1\) in the interior or the cycle \(C'_{12}\) formed by \(P'_1 \cup P_2\) contains \(P_1\) in the interior. Suppose the former holds and let \(S_{a}\) be any occupied square in the interior of \(E_{r_1}.\) We know that \(E_{r_1} = D_a\) is the outermost boundary cycle containing \(S_a\) and satisfies conditions (a), (b) and (c) mentioned in Lemma~\ref{outer}. The cycle \(C_{12}\) also contains \(S_{a}\) in the interior and thus satisfies condition (a). Moreover, it contains at least one edge in the exterior of \(E_{r_1}\) contradicting the fact that \(E_{r_1}\) satisfies (c). Thus \(H_{cyc}\) is acyclic. \(\qed\)

\section{Proofs of Theorem~\ref{thm2} and Theorem~\ref{thm1}}\label{pf1}
\emph{Proof of Theorem~\ref{thm2}}: Let \(D_0\) be the outermost boundary cycle containing the square \(S_0\) as in Lemma~\ref{outer}. It satisfies the conditions (i) and (ii) in the statement of the theorem and is unique and thus \(C^+_{out} = D_0.\) \(\qed\)

\emph{Proof of Theorem~\ref{thm1}}: If every edge of \(C_1\) is either on \(C_2\) or contained in the interior of \(C_2,\) then the desired cycle \(C_3 = C_2.\) If similarly, \(C_2\) is completely contained in \(C_1,\) we set \(C_3 = C_1.\) So we suppose that \(C_1\) contains at least one edge in the exterior of \(C_2\) and \(C_2\) also contains at least one edge in the exterior of \(C_1.\)

We start with cycle \(C_1\) and in the first step, identify a path of \(C_2\) contained in the exterior of \(C_1.\) Set \(C_{1,0} := C_1 = (u_0,u_1,...,u_{t-1},u_0)\) and  \(C_2 = (v_0,v_1,...,v_{m-1},v_0).\) For later notation, we define \(u_{k} = u_{k \mod t}\) if \(k \leq 0\) or \(k \geq t\) and \(v_{k} = v_{k \mod m}\) if \(k \leq 0\) or \(k \geq m.\)

Start from some vertex of \(C_{1,0},\) say \(u_0,\) and look for the first intersection point that contains an exterior edge of \(C_2;\) i.e., an edge of \(C_2\) that lies in the exterior of \(C_1.\) Let \[j_1 = \min\{j \geq 0 : u_j \in C_{1,0} \text{ and } u_j \text{ is an endvertex of an exterior edge of } C_2 \}\] and let \(v_{i_1} = u_{j_1}.\) We suppose that the edge of \(C_2\) with endvertices \(v_{i_1}\) and \(v_{i_1+1}\) lies in the exterior of \(C_1.\) Let \(r_1 = \min\{i \geq i_1+1 : v_i \in C_{1,0}\}\) be the next time the cycles meet and define \(P_1 = (v_{i_1},v_{i_1+1},...,v_{r_1}).\)

We note that none of the vertices \(v_{j}, i_1 +1 \leq j \leq r_1-1\) belong to \(C_{1,0}.\) If \(v_{i_1} = v_{r_1},\) then \(P_1\) is a cycle containing the edges of \(C_2\) and thus \(P_1 = C_2.\) Since \(C_1\) and \(C_2\) contain more than one vertex in common, this cannot happen. Thus \(P_1\) is a path and all edges of \(P_1\) are in the exterior of~\(C_{1,0}.\)

We then construct an outermost cycle from \(C_{1,0}\) and \(P_1\) as follows. Split \(C_{1,0}\) into two segments based on intersection with \(P_1.\) Suppose \(P_1\) meets \(C_{1,0}\) at \(u_{a_1}\) and \(u_{b_1}.\) We let \(C'_{1,0} = (u_{a_1},u_{a_1+1},...,u_{b_1})\) and \(C^{''}_{1,0} = (u_{a_1}, u_{a_1-1},...,u_{b_1}).\)  If the interior of \(C'_{1,0} \cup P_1\) contains the interior of \(C''_{1,0} \cup P_1\) as in Figure~\ref{fig_first_case}, we set \(C_{1,1} = C'_{1,0} \cup P_1\) to be the cycle obtained in the first iteration by the concatenation of the paths \(C'_{1,0}\) and \(P_1.\) Here \(C''_{1,0}\) is the segment \(ADC,\) the path \(C'_{1,0}\) is the segment \(ABC\) and the path \(P_1\) is denoted \(AEC.\) Else necessarily we have that the interior of \(C^{''}_{1,0} \cup P_1\) contains the interior of \(C'_{1,0} \cup P_1\) and we set \(C_{1,1} = C^{''}_{1,0} \cup P_1.\) Since \(P_1 \neq \emptyset,\) we have that \(C_{1,1}\) contains at least one exterior edge.

We then perform the same procedure as above on the cycle \(C_{1,1}\) and continue this process for a finite number of steps to obtain the final cycle~\(C_{1,n}.\) For each \(j, 1 \leq j \leq n,\) we have that the cycle \(C_{1,j}\) satisfies the following properties:\\
(1) the cycle \(C_{1,j}\) contains only edges from \(C_1\) and \(C_2,\)\\
(2) every edge of \(C_1\) either belongs to \(C_{1,j}\) or is contained in the interior of~\(C_{1,j},\) \\
(3) the cycle \(C_{1,j}\) contains at least one exterior edge of \(C_2\) and \\
(4) the interior of \(C_1\) is contained in \(C_{1,j}.\) \\
In particular, the above properties hold true for the final cycle \(C_{1,n}.\) If there exists an edge \(e\) of \(C_2\) in the exterior of \(C_{1,n},\) then the edge \(e\) belongs to a path \(P_e\) of \(C_2\) containing edges exterior to \(C_1.\) The path \(P_e\) must meet \(C_1\) and thus there exists an edge of \(C_2\) that lies in the exterior of \(C_1\) and contains an endvertex of \(C_1.\) But then the above procedure would not have terminated and thus we also have:\\
(5) every edge of \(C_2\) either belongs to \(C_{1,n}\) or is contained in the interior of~\(C_{1,n}.\)

Thus property (ii) stated in the result holds true and we need to see that (i) holds. For that we first prove uniqueness of the cycle \(C_{1,n}\) obtained above. Suppose there exists another cycle \(D'\) satisfying properties (1), (2) and (5) above. If \(D'\) contains an edge \(e'\) (which must necessarily belong to \(C_1\) or \(C_2\)) in the exterior of \(C_{1,n},\) it contradicts the fact that \(C_{1,n}\) satisfies (2) and (5). If \(D'\) is completely contained in the interior of \(C_{1,n}\) and is not equal to \(C_{1,n},\) then there is at least one edge of \(C_{1,n}\) (which belongs to \(C_1\) or \(C_2\)) that lies in the exterior of~\(D',\) contradicting the assumption that \(D'\) satisfies~(2)~and~(5).

Thus any cycle satisfying properties (1), (2) and (5) is unique. We recall that \(C_1\) also contains an edge in the exterior of \(C_2.\) Suppose now we start from \(C_{2,0} := C_2\) and identify segments of \(C_1\) lying in the exterior of \(C_2\) and perform the same iterative procedure as above to obtain a final cycle \(C_{2,m}.\) This cycle must also satisfy (1), (2) and (5) and hence \(C_{2,m} = C_{1,n}.\) Moreover, \(C_{2,m}\) satisfies:\\
(3\('\)) the cycle \(C_{2,m}\) contains at least one exterior edge of \(C_1\) and \\
(4\('\)) the interior of \(C_2\) is contained in \(C_{2,m}.\)

Thus the cycle \(C_{1,n}\) is unique and satisfies properties (i) and (ii) stated in the result. \(\qed\)

\renewcommand{\theequation}{\thesection.\arabic{equation}}

\subsection*{Acknowledgement}
I thank Professor Rahul Roy for crucial comments and NISER for my fellowship.

\bibliographystyle{plain}

\end{document}